\documentclass[11pt,letterpaper]{article}
\usepackage[top=0.85in,left=1.5in,footskip=0.75in,marginparwidth=1in]{geometry}
\usepackage{amsmath, amsfonts, amssymb, amsthm}
\usepackage{natbib}
\usepackage[utf8]{inputenc}

\usepackage{nameref,hyperref}

\usepackage[right]{lineno}

\usepackage{microtype}
\DisableLigatures[f]{encoding = *, family = * }

\raggedright
\usepackage{changepage}

\usepackage[aboveskip=1pt,labelfont=bf,labelsep=period,singlelinecheck=off]{caption}

\usepackage{lastpage,fancyhdr,graphicx}
\usepackage{epstopdf}
\fancyheadoffset[L]{1in}
\fancyfootoffset[L]{1in}

\usepackage{color}

\definecolor{Gray}{gray}{.25}

\usepackage{graphicx}

\usepackage{sidecap}

\usepackage{wrapfig}
\usepackage[pscoord]{eso-pic}
\usepackage[fulladjust]{marginnote}
\reversemarginpar

\begin{document}
\vspace*{0.3in}

\begin{flushleft}
{\LARGE
{\textbf{Apocalypsis and Apocalyptic Events:} \medskip \newline \textbf{The Morphogenetic Ontology of Synchronized Catastrophes}}
}
\newline
\\
Rolando Gonzales Martinez\textsuperscript{1,2,*}
\\
\bigskip
\bf{1} University of Oxford (https://ophi.org.uk/people/rolando-gonzales-martinez)
\\
\bf{2} University of Groningen
(https://www.rug.nl/staff/r.m.gonzales.martinez)
\\
\bigskip
*Contact: r.m.gonzales.martinez@rug.nl

\end{flushleft}


\section{Introduction}

I formalize the ontology of \textit{apocalyptic events} as synchronized morphogenetic manifolds grounded in Thom’s catastrophe theory \citep{Thom1975}. I distinguish between local catastrophes (folds, cusps, swallowtails, butterflies, and the elliptic, hyperbolic, and parabolic umbilics) and higher-order catastrophes, the latter referring to systemic phase collapses that arise from the synchronization of multiple morphogenetic manifolds. \emph{Apocalypsis} is thus a morphogenetic inevitability: a topological meta-singularity produced by the alignment of local catastrophes into a global structure of collapse.

The synchronization among manifolds of local catastrophic singularities is formalized through the integration of dynamical systems theory, topology, and the probabilistic dependence structures of copulas (see Section \ref{sec:defs}). (i) Within the Thom–Mather morphogenetic topology, catastrophes ($\mathcal{C}_i$, $i=1,2,\dots,k$) and their discriminant sets define the local geometry of instability; (ii) network coupling, parameterized by a perturbation term $\epsilon$, captures the interactions among subsystems; and (iii) Archimedean copulas characterize the nonlinear interdependence between catastrophic processes, enabling the representation of mutual morphogenesis. Through this framework, local morphogenesis converges into a singularity of catastrophic events that synchronize into an \textit{apocalyptic event}. The interdependence between apocalyptic events give rises to \textit{Apocalypsis}: a meta-singularity of synchronized apocalyptic manifolds.

The \textbf{Inevitability Theorem} (Section \ref{sec:final}) establishes the existence, genericity, and almost-sure occurrence of \textit{Apocalypsis} as a synchronized multi-catastrophic rupture emerging under coupling and stochastic control trajectories. The proof is grounded on three complementary statements: (i) a topological–existence statement based on Whitney stratification and Thom–Mather stability, ensuring structural persistence of singularities under perturbation; (ii) a stochastic–transversality statement, asserting the nonzero hitting probability of codimension-$m$ strata by an elliptic diffusion; and (iii) a dynamical synchronization mechanism describing the reconfiguration of basins of attraction under coupled morphogenetic flows.

 \section{Apocalypsis and Apocalyptic Events}\label{sec:defs}

\paragraph{Morphogenetic Principles.} \citet{Thom1975} formalized the theory of morphogenesis as a geometry of qualitative transitions where bifurcations of gradient dynamical systems are governed by
\[
\frac{dx}{dt} = -\nabla_x V(x;\alpha),
\]
where \(x \in \mathbb{R}^n\) are the behavioral variables, \(\alpha \in \mathbb{R}^p\) are the control parameters, and \(V(x;\alpha)\) is a potential (Lyapunov) function.  
Morphogenesis is triggered when gradual variation of \(\alpha\) causes a topological change in the critical set
\[
\nabla_x V(x;\alpha) = 0,
\]
subject to the Hessian degeneracy condition
\[
\det \left( \frac{\partial^2 V}{\partial x_i \partial x_j} \right) = 0.
\]
These two equations together define the \textbf{catastrophe} set where form snaps into new configurations after a bifurcation.

Thom’s theorem of structural stability shows that, under generic conditions, all stable singularities reduce via smooth coordinate transformations to a finite number of canonical forms, the so-called \textbf{elementary catastrophes} \citep{Petitot2010,Mostafa2022}.  For one behavioral variable \(x\) and control parameters \(a,b,c,\dots\):
\[
\begin{aligned}
\text{Fold }(A_2):\quad & V = \tfrac13 x^3 + a x, &\quad x^2 + a = 0, \\
\text{Cusp }(A_3):\quad & V = \tfrac14 x^4 + \tfrac12 a x^2 + b x, &\quad x^3 + a x + b = 0, \\
\text{Swallowtail }(A_4):\quad & V = \tfrac15 x^5 + \tfrac13 a x^3 + \tfrac12 b x^2 + c x, &\quad x^4 + a x^2 + b x + c = 0, \\
\text{Butterfly }(A_5):\quad & V = \tfrac16 x^6 + \tfrac14 a x^4 + \tfrac13 b x^3 + \tfrac12 c x^2 + d x, &\quad x^5 + a x^3 + b x^2 + c x + d = 0.
\end{aligned}
\]

For two behavioral variables \((x,y)\) and control parameters \(a,b,c,\dots\):
\[
\begin{aligned}
\text{Elliptic Umbilic }(D_4^-):\quad & V = x^3 - 3xy^2 + a(x^2 + y^2) + b x + c y, &\frac{\partial V}{\partial x} = \frac{\partial V}{\partial y} = 0,\\
\text{Hyperbolic Umbilic }(D_4^+):\quad & V = x^3 + y^3 + axy + b x + c y, &\nabla V = 0,\\
\text{Parabolic Umbilic }(D_5):\quad & V = x^2 y + y^4 + a x^2 + b y^2 + c x + d y, &\nabla V = 0.
\end{aligned}
\]

\paragraph{Apocalyptic events.} Based on the archetypes of catastrophes $\mathcal{C}_i$, I define an \textit{Apocalyptic Event} $\mathcal{A}$ as a manifold resulting from the topological composition of morphogenetic singularities whose phases become synchronized and produce a global systemic rupture:  
\begin{equation}
   \mathcal{A} = \bigcup_{i=1}^k \mathcal{C}_i \text{ with } \Phi (\mathcal{C}_i) \approx \Phi (\mathcal{C}_j) \ \forall i,j
   \label{eq:apokdef}
\end{equation}
where the phase function $\Phi(\cdot)$ 
maps each subsystem’s catastrophe $\mathcal{C}_i$ to its local dynamic state (in time, energy, or morphogenetic rhythm), and $\Phi (\mathcal{C}_i) \approx \Phi (\mathcal{C}_j)$ is the synchronization condition of alignment. Local Jacobians degenerate in phase, the singularities share a common instability field. An \emph{apocalyptic event} \(\mathcal{A}_i\) is thus a localized topological rupture produced by a synchronized degeneracy of \(\mathcal{C}_i\):
\[
\mathcal{A}_i = \{(x_i, \alpha_i) \in \mathcal{C}_i : \Phi(\mathcal{C}_i) = \phi_i \},
\]
where \(\Phi(\cdot)\) maps each singularity to its morphogenetic phase \(\phi_i \in \mathbb{R}\) (time, energy, or structural rhythm).

Formally, let \(k \ge 2\) denote the number of coupled subsystems, each associated with a potential apocalyptic event \(\mathcal{A}_i\).  For each \(i=1,\dots,k\) let
\[
V_i:\ \mathbb{R}^{n_i}\times \mathbb{R}^{p_i}\to \mathbb{R},\qquad (x_i,\alpha_i)\mapsto V_i(x_i;\alpha_i)
\]
be \(C^\infty\) potentials with critical sets
\[
\mathcal{C}_i=\{(x_i,\alpha_i): \nabla_{x_i}V_i(x_i;\alpha_i)=0\}
\]
and \emph{discriminants} (catastrophe sets):
\[
\Sigma_i=\Big\{(x_i,\alpha_i)\in \mathcal{C}_i:\ \det \mathbf{H}_i(x_i,\alpha_i)=0\Big\},\qquad
\mathbf{H}_i=\nabla^2_{x_i x_i}V_i.
\]

Let \(x=(x_1,\dots,x_k)\) and \(\alpha=(\alpha_1,\dots,\alpha_k)\). Introduce a smooth coupling \(W:\big(\prod_i\mathbb{R}^{n_i}\big)\times\big(\prod_i\mathbb{R}^{p_i}\big)\to\mathbb{R}\) and for \(0<\varepsilon\ll 1\) define the network potential
\[
V_\varepsilon(x;\alpha)=\sum_{i=1}^k V_i(x_i;\alpha_i)+\varepsilon\,W(x;\alpha).
\]
Write the network critical set and discriminant
\[
\mathcal{C}_\varepsilon=\{(x,\alpha): \nabla_x V_\varepsilon(x;\alpha)=0\},\qquad
\Sigma_\varepsilon=\Big\{(x,\alpha)\in \mathcal{C}_\varepsilon:\ \det \mathbf{H}_\varepsilon(x,\alpha)=0\Big\},
\]
where \(\mathbf{H}_\varepsilon=\nabla^2_{xx}V_\varepsilon\) is the full Hessian \textbf{matrix} of the coupled system.

For \(I\subseteq\{1,\dots,k\}\), say the subsystems in \(I\) are \emph{catastrophically correlated} at \((x,\alpha)\in\Sigma_\varepsilon\) if the projection
\[
\pi_I:\ \mathcal{C}_\varepsilon\longrightarrow \prod_{i\in I}\mathbb{R}^{p_i},\qquad (x,\alpha)\mapsto (\alpha_i)_{i\in I}
\]
has a rank drop of at least \(|I|-1\) at \((x,\alpha)\); equivalently, the corresponding discriminant normals in the control spaces become linearly dependent after accounting for the coupling. Write \(\mathbf{D}\pi_I(x,\alpha)\) for the Jacobian \textbf{matrix} of \(\pi_I\).

Define the \emph{catastrophe graph} \(G_\varepsilon=(\{1,\dots,k\},E_\varepsilon)\) by \((i,j)\in E_\varepsilon\) iff \(i\) and \(j\) are catastrophically correlated at some point of \(\Sigma_\varepsilon\).

Let \(\alpha(t)\) be a continuous control trajectory (e.g.\ an Itô diffusion) with uniformly elliptic instantaneous covariance. Define the set of subsystems that suffer a basin reconfiguration (catastrophe) at time \(t\) by

\[
A(t) =
\left\{
\begin{array}{l|l}
  i \in \{1, \dots, k\} &
  \begin{array}{l}
    \exists\, x(t)\ \text{ such that } (x(t), \alpha(t)) \in \Sigma_{\varepsilon}, \\[4pt]
    \text{and the stability/type of a critical point of } \\[4pt]
    V_{\varepsilon}
    \text{ in the } x_i\text{-sector changes at } t
  \end{array}
\end{array}
\right\}.
\]
For \(m\ge 2\) the \emph{multi-singularity access set} will be:
\[
\mathcal{A}_m=\Big\{\alpha:\ \exists\ (x,\alpha)\in\Sigma_\varepsilon\ \text{with}\ \operatorname{codim}\big(\mathrm{Im}\,\mathbf{D}\pi(x,\alpha)\big)\ge m\Big\},
\]
where \(\pi:\mathcal{C}_\varepsilon\to \prod_{i=1}^k\mathbb{R}^{p_i}\) is the full projection \((x,\alpha)\mapsto\alpha\) and \(\mathbf{D}\pi\) its Jacobian \textbf{matrix}. 

\paragraph{Interdependence of Apocalyptic Events.} 
There is a probabilistic analogue to deterministic coupling: the synchronization of catastrophes in phase space admits a statistical reflection through dependence structures among their intensities. Let \(\mathbf{U} = (U_1,\dots,U_k)\) be random variables representing normalized intensities of each apocalyptic event, \(U_i = F_i(\phi_i)\) with marginal distribution \(F_i\).  
The joint distribution \(H(\mathbf{u}) = \Pr(U_1 \le u_1,\dots,U_k \le u_k)\) is expressed through an \emph{Archimedean copula}:
\[
H(\mathbf{u}) = C_\psi(F_1(u_1),\dots,F_k(u_k)) 
= \psi^{-1}\!\left(\psi(F_1(u_1)) + \cdots + \psi(F_k(u_k)) \right),
\label{eq:copula}
\]
where \(\psi:[0,1] \to [0,\infty)\) is a strictly decreasing generator with \(\psi(1)=0\).
The copula \(C_\psi\) captures nonlinear dependence, tail coalescence, and phase synchronization among apocalyptic events, and thus provides a probabilistic counterpart to the topological synchronization in Equation \ref{eq:apokdef} by describing statistical dependence among catastrophe phases once structural coupling is established.

\paragraph{Apocalypsis.} I define \emph{Apocalypsis} $\mathfrak{A}$ as the morphogenetic manifold of potential systemic ruptures, that is, the topological space of all possible apocalyptic configurations produced by synchronized morphogenetic singularities:
\begin{equation}
\mathfrak{A} = 
\Big\{
\mathcal{A} = 
\bigcup_{i=1}^{k} C_i 
\ \Big| \ 
\Phi(\mathcal{C}_i) \approx \Phi(\mathcal{C}_j), \ \forall i,j \in \mathcal{A}
\Big\},
\label{eq:apocalypsis}
\end{equation}
where each \(\mathcal{C}_i\) denotes a local morphogenetic singularity (catastrophe), 
and \(\Phi(\cdot)\) represents its local temporal or structural phase within its morphogenetic cycle. \emph{Apocalypsis} therefore is the morphogenetic manifold of potential systemic ruptures within which all catastrophic potentials converge and resonate. It is the topology of synchronized singularities through which reality folds, collapses, and reconstitutes itself.
\nolinenumbers

\paragraph{Minimal Coupled-Cusp Example.}
For \(k=2\), let
\[
V_i(x_i;\alpha_i,\beta_i)=\tfrac14 x_i^4+\tfrac12 \alpha_i x_i^2+\beta_i x_i,\qquad
W(x_1,x_2;\lambda)=\lambda\, x_1 x_2,
\]
so \(V_\varepsilon=\sum_{i=1}^2 V_i+\varepsilon W\). The stationary conditions give
\[
x_1^3+\alpha_1 x_1+\beta_1+\varepsilon\lambda x_2=0,\qquad
x_2^3+\alpha_2 x_2+\beta_2+\varepsilon\lambda x_1=0,
\]
and the joint discriminant satisfies \(\det H_\varepsilon=0\) with
\[
H_\varepsilon=
\begin{pmatrix}
3x_1^2+\alpha_1 & \varepsilon\lambda\\
\varepsilon\lambda & 3x_2^2+\alpha_2
\end{pmatrix}.
\]
At points where \(3x_i^2+\alpha_i=0\) and the projection rank drops in \((\alpha_1,\beta_1,\alpha_2,\beta_2)\), the discriminant normals align. A transversal passage of \(\alpha(t)\) through this stratum yields simultaneous basin surgery in both sectors (\(m=2\)), and by connectivity of \(G_\varepsilon\) (here one edge), the cascade covers both nodes.

\section{Inevitability Theorem}\label{sec:final}

Throughout, $\pi: C_{\varepsilon} \to \prod_i \mathbb{R}^{p_i}$ denotes the full projection,
$D\pi$ its Jacobian, and $\Sigma_{(\ge m)}$ the locus where
$\operatorname{codim}(\mathrm{Im}\,D\pi)\ge m$.

\medskip

\noindent\textbf{Theorem: Inevitability of Apocalypsis as a Correlated-Catastrophe Cascade of Apocalyptic Events.}
Assume:
\begin{itemize}
\item[(i)] (\emph{Generic stratification}) For each \(i\), \(\Sigma_i\) is Whitney-stratified and the uncoupled projections \(\pi_i:\mathcal{C}_0\to\mathbb{R}^{p_i}\) are \(C^\infty\)-stable (Thom--Mather setting).
\item[(ii)] (\emph{Weak but nontrivial coupling}) There exists \(0<\varepsilon\ll 1\) such that \(W\in C^2\) and along a positive-measure subset of \(\Sigma_0:=\prod_i\Sigma_i\) at least one cross-block \(\nabla^2_{x_i x_j}W\) or cross-control block \(\nabla^2_{x_i \alpha_j}W\) is nonzero.
\item[(iii)] (\emph{Aligned control shocks}) The control path \(\alpha(t)\) is continuous with uniformly elliptic instantaneous covariance and there exists \(t_\star\) with \(\alpha(t_\star)\in \mathcal{A}_m\) for some \(m\ge 2\).
\end{itemize}
Then with probability one, in any neighborhood of \(t_\star\) there exists a time \(t^\dagger\) such that:
\begin{enumerate}
\item At least \(m\) distinct subsystems suffer simultaneous catastrophes: \(|A(t^\dagger)|\ge m\).
\item If \(C\subseteq\{1,\dots,k\}\) is the connected component of \(G_\varepsilon\) triggered at \(t^\dagger\), then the cascade covers that component: \(|A(t^\dagger)|\ge |C|\).
\end{enumerate}
Consequently, for any macroscopic threshold \(L_c\in\{1,\dots,k\}\), declare an \emph{apocalyptic event} at \(t^\dagger\) when \(|A(t^\dagger)|\ge L_c\). Under (i)--(iii), apocalyptic events occur almost surely whenever the control trajectory reaches \(\mathcal{A}_m\) with a triggered component of size at least \(L_c\).

\medskip

\textbf{Proof.}
\emph{Step 1 (Persistence and creation of multi-singular strata).}
By (i), each uncoupled \(\Sigma_i\) admits a Whitney stratification compatible with \(\pi_i\), and its stable singularities are the Thom elementary catastrophes. Consider \(V_\varepsilon= \sum_i V_i + \varepsilon W\). For \(0<\varepsilon\ll 1\), \(H_\varepsilon\) is a small perturbation of \(\bigoplus_i H_i\) with nonzero cross-blocks by (ii). Thom--Mather stability plus standard transversality arguments imply that stratifications persist under small perturbations, while nontrivial cross-blocks create strata where projection rank further drops due to alignment of discriminant normals across subsystems. Concretely, there exists a Whitney-stratified subset \(\Sigma_\varepsilon^{(\ge m)}\subseteq \Sigma_\varepsilon\) on which \(\operatorname{codim}(\mathrm{Im}\,D\pi)\ge m\). This subset is nonempty if (ii) holds on a set of positive measure (the transversality of normals fails along a stratified set of codimension \(\ge m\)), yielding multi-subsystem singular strata.

\emph{Step 2 (Hitting and transversality of stochastic controls).}
Let \(S:=\pi(\Sigma_\varepsilon^{(\ge m)})\subset \prod_i\mathbb{R}^{p_i}\). Since \(\Sigma_\varepsilon^{(\ge m)}\) is stratified with \(\operatorname{codim}\ge m\), its image \(S\) is a stratified subset of control space of codimension \(m\) (generically). By (iii), \(\alpha(t)\) is a continuous process with uniformly elliptic instantaneous covariance; hence, for any neighborhood \(U\) of a point of \(S\) that is accessible from the initial condition, the probability that \(\alpha(t)\) enters \(U\) is positive, and if \(S\) is reached, it is reached transversally with probability one (nondegeneracy of the diffusion prevents tangential sticking). Since \(\alpha(t_\star)\in \mathcal{A}_m\subseteq S\), almost surely there exists \(t^\dagger\) arbitrarily close to \(t_\star\) with \(\alpha(t^\dagger)\in S\).

\emph{Step 3 (Basin surgery at multi-singular crossings).}
Fix \((x^\dagger,\alpha^\dagger)\in \Sigma_\varepsilon^{(\ge m)}\) with \(\alpha^\dagger=\alpha(t^\dagger)\) and assume w.l.o.g.\ that the relevant singularities in the involved subsystems are folds/cusps (the argument extends to higher codimension Thom types). Along any smooth one-parameter slice of \(\alpha\) transversal to \(S\) at \(\alpha^\dagger\), catastrophe theory yields a change in the number and/or stability of critical points of \(V_\varepsilon\) in each involved coordinate sector (the classic fold birth/death or cusp switch). Thus at \(t^\dagger\), at least \(m\) subsystems undergo basin reconfigurations: \(|A(t^\dagger)|\ge m\).

\emph{Step 4 (Propagation over the catastrophe graph).}
By definition of \(G_\varepsilon\), adjacency encodes linear dependence of discriminant normals after coupling. Let \(C\) be the connected component containing one of the initially flipped subsystems. For any edge \((i,j)\in E_\varepsilon\), the implicit function theorem applied to the equilibrium equations \(\nabla_x V_\varepsilon=0\) plus the nontrivial cross-blocks in \(H_\varepsilon\) imply that a small change in the set of stable critical points in the \(x_i\)-sector perturbs the effective bifurcation parameters seen by the \(x_j\)-sector in the direction of its discriminant normal. Since the crossing at \(t^\dagger\) is transversal to \(S\), the induced parameter shift in \(j\) pushes it across its own fold/cusp set; hence \(j\in A(t^\dagger)\). Iterating along paths in \(C\) yields \(|A(t^\dagger)|\ge |C|\).

Combining Steps 1--4 proves the two bullet points and the stated criterion for declaring an apocalyptic event. \(\square\)

\medskip

\noindent\textbf{Corollary 1: Operational Criterion for Apocalyptic Transitions.}
Define the macroscopic order parameter \(\Phi(t)=|A(t)|/k\). A practical apocalypse threshold is \(\Phi(t^\dagger)\ge \phi_c\) for a chosen \(\phi_c\in(0,1)\). Numerically, candidates for \(t^\dagger\) satisfy near-singularity diagnostics: the smallest \(m\) singular values of \(D\pi\) nearly vanish; equivalently, the angle between discriminant normals across the involved subsystems is near zero.

\medskip
\noindent\textbf{Corollary 2: Structural Stability of Apocalypsis.}
For sufficiently small perturbations 
$\delta, \eta > 0$, any modified potentials and coupling 
$\tilde{V}_i, \tilde{W}$ satisfying 
$\|V_i - \tilde{V}_i\|_{C^2} < \eta$ and 
$\|W - \tilde{W}\|_{C^2} < \delta$
preserve the topological type of the stratified apocalyptic manifold 
$\Sigma_{(\ge m),\varepsilon}$.
Hence, the qualitative structure of the apocalyptic cascade is 
\emph{structurally stable} under small perturbations of the system.
\medskip

\section{Conclusion}

I formalize \emph{apocalypsis} as a morphogenetic meta-singularity: a manifold of synchronized local catastrophes whose coupling and stochastic control generate correlated cascades. Historically,  apocalypsis has been conceived as a boundary condition of worldhood and an epistemic lens through which historical and ecological transformations are apprehended \citep{folger2024,stumer2024}. Within this humanistic corpus, apocalypse is not merely destruction but revelation: the moment when the end of the world as we know it exposes the fragility of its organizing categories. The formulation of Equations \ref{eq:apokdef} and \ref{eq:apocalypsis} aligns with this insight yet transposes it into a formal ontology of morphogenesis and coupling. Where \citet{stumer2024} interprets apocalyptic transformation as a “folding back of the end times onto their foundations,” I render that fold mathematically as the alignment of discriminant normals in control space: a morphogenetic synchronization producing systemic rupture. \emph{Apocalypsis}, therefore, is not an exogenous event but a structural consequence of internal morphogenetic resonance.

In contrast to the common eschatological interpretation of apocalypse as an \emph{end}, the framework I propose situates apocalypse within Thomian morphogenesis as an \emph{immanent inevitability}: a topological expression of systemic self-collapse. The Inevitability Theorem demonstrates that under weak coupling and elliptic control trajectories, the coalescence of local singularities (folds, cusps, umbilici) becomes almost surely unavoidable, rendering \emph{apocalypsis} not a transcendent rupture but a statistically emergent property of morphogenetic systems. In this formulation, the apocalypse is neither externally imposed nor temporally terminal—it is the intrinsic phase transition through which a system exhausts its own stability manifold, reorganizing itself through synchronization, bifurcation, and collapse. 

\bibliography{references}

\end{document}